\documentclass{amsart}
\usepackage{amsmath, amssymb,verbatim}
\newtheorem{theorem}{Theorem}

\newtheorem{lemma}{Lemma}

\author[BOUHALI]{Aissa Bouhali}
\address{Département de Mathématiques, Ecole Normale Supérieure de Laghouat, Algeria.
    Laboratory of Pure and Applied Mathematics, University of Laghouat, Algeria.}
\email{aissa.bouhali@ens-lagh.dz}
\author[LOUHICHI]{Issam Louhichi }
\address{Department of Mathematics \& Statistics, College of Arts \& Sciences, American University of Sharjah, P.O.Box 2666, Sharjah, UAE.}
\email{ilouhichi@aus.edu}

\begin{document}
\title[Commutants of the sum of two Toeplitz operators ]
{Commutants of the sum of two quasihomogeneous Toeplitz operators}

\date{\today} 
\begin{abstract}
	A major open question in the theory of Toeplitz operator on the Bergman space of the unit disk of the complex plane is the complete characterization of the set of all Toeplitz operators that commute with a given operator. In \cite{al}, the authors showed that when a sum $S=T_{e^{im\theta}f}+ T_{e^{il\theta}g}$, where $f$ and $g$ are radial functions,  commutes with a sum $T=T_{e^{ip\theta}r^{(2M+1)p}}+ T_{e^{is\theta}r^{(2N+1)s}}$, then $S$ must be of the form $S=cT$, where $c$ is a constant. In this article, we will replace $r^{(2M+1)p}$ and $r^{(2N+1)s}$ with $r^n$ and $r^d$, where $n$ and $d$ are in $\mathbb{N}$, and we will show that the same result holds.
\end{abstract}	
\maketitle 
\begin{center}
	\textit{Dedicated to Rao V. Nagisetty (1938-2024).}
\end{center}
\section{Introduction}
In the complex plane $\mathbb{C}$, let $\mathbb{D}$ be the open unit disk, and $dA=rdr\frac{d\theta}{\pi}$ be the normalized Lebesgue measure on $\mathbb{D}$, where $(r,\theta)$ are the polar coordinates. The Hilbert space $L^2(\mathbb{D},dA)$ consists  of all square-integrable functions on $\mathbb{D}$ with respect to the measure $dA$. 

The classical unweighted Bergman space $L^2_a(\mathbb{D})$ is the closed subspace of $L^2(\mathbb{D},dA)$ consisting of all analytic functions on $D$. Moreover, the set $\{z^n: n=0,1,2,\ldots\}$ is an orthogonal basis for $L^2_a(\mathbb{D})$. Since $L^2_a(\mathbb{D})$ is closed, the orthogonal projection $P$ from $L^2(\mathbb{D},dA)$ onto $L^2_a(\mathbb{D})$, known as the Bergman projection, is well-defined.  For more information on the theory of Bergman spaces, see \cite{h} .

For a function $f\in L^2(\mathbb{D},dA)$, we define the Toeplitz operator $T_f$ from $L^2_a(\mathbb{D})$ to itself by $T_f(g)=P(fg)$ whenever $fg$ is in $L^2(\mathbb{D}, dA)$. The function $f$ is called the symbol of the Toeplitz operator $T_f$. From this definition, it is clear that bounded analytic functions are within the domain of $T_f$, making $T_f$ densely defined on $L^2_a(\mathbb{D})$. Furthermore,  if the symbol $f$ is bounded on $\mathbb{D}$, then $T_f$ is bounded and $||T_f||\leq ||f||_{\infty}$.

Over the past fifty years, various algebraic properties of Toeplitz operators have been extensively studied. Nevertheless, the problem of describing the commutant of a given Toeplitz operator, i.e., the set of all Toeplitz operators that commute with it in the sense of composition, remains largely unresolved. Very little is known about when $T_fT_g=T_gT_f$ for "general symbols" $f$ and $g$. For more information on the commutativity of Toeplitz operators on $L^2_a(D)$ and the results obtained so far for certain specific classes of symbols, see \cite{sl, al, acr, cr, lt, l, lr, lry, lz, rv, sz,y}.

In this work, we address the class of  so-called quasihomogeneous Toeplitz operators. A symbol $f$ is said to be quasihomogeneous of degree  $p$ (an integer) if $f(re^{i\theta})=e^{ip\theta}\phi(r)$, where $\phi$ is a radial function. In this case, the associated Toeplitz operator $T_f$ is also referred to as a quasihomogeneous Toeplitz operator of degree $p$ (see \cite{cr, lr, lsz}). The motivation for studying this family of symbols is that $L^2(\mathbb{D},dA)$ can be expressed as $L^2(\mathbb{D},dA)=\bigoplus_{k\in\mathbb{Z}}e^{ik\theta}\mathcal{R}$, where $\mathcal{R}$ is the space of square-integrable radial functions on $[0,1)$ with respect to the measure $rdr$. Thus every function $f\in L^2(\mathbb{D},dA)$ has the polar decomposition $f(z)=f(re^{i\theta})=\sum_{k\in\mathbb{Z}}e^{ik\theta}f_k(r)$, where the $f_k$'s are radial functions. We believe that studying quasihomogeneous Toeplitz operators will help us  characterize the commutant of Toeplitz operators with more general symbols.

This present work is motivated by \cite[Theorem 3.1,~p.52]{al}. For our purposes, we will state and summarize this theorem as follows:
\begin{theorem}\label{thm4}
	Let $\phi(r)=r^{(2M+1)p}$ and $\psi(r)=r^{(2N+1)s}$, where $p<s$ and $M$ and $N$ are integers greater or equal to $1$. Suppose there exist $m, l\in\mathbb{N}$ and nontrivial radial functions $f$ and $ g$ such that  the following two hypotheses are satisfied
	\begin{displaymath}
		\left\{\begin{array}{ll}
			(H1)& T_{e^{im\theta}f}+ T_{e^{il\theta}g}
			\textrm{ commutes with } T_{e^{ip\theta}\phi}+T_{e^{is\theta}\psi},\\
			(H2)& 1\leq p<s, \ 1\leq m<l, \textrm{ and } l+p=m+s.
		\end{array}\right.
	\end{displaymath}
	Then $m=p,\ l=s$ and
	$$T_{e^{im\theta}f}+T_{e^{il\theta}g}=c\left(T_{e^{ip\theta}\phi}+T_{e^{is\theta}\psi}\right)$$ for some constant $c$.
\end{theorem}
Our goal here is to replace the radial components of the symbols $e^{ip\theta}\phi$ and $e^{is\theta}\psi$ in Theorem 1, specifically, $\phi(r)=r^{(2M+1)p}$ and $\psi(r)=r^{(2N+1)s}$, with $\phi(r)=r^n$ and $\psi(r)=r^d$, where $n$ and $d$ are in $\mathbb{N}$. In other words, we aim to eliminate the condition that the power of $r$ in $\phi$ and $\psi$ is an odd number multiple of the quasihomogeneous degrees $p$ and $s$, respectively. Despite this modification, we will still be able to apply the result in \cite{lrp} regarding the existence of roots for Toeplitz operators. This result states that there exist radial functions $\widetilde{\phi}$ and $\widetilde{\psi}$ such that $T_{e^{ip\theta}r^n}=\left(T_{e^{i\theta}\widetilde{\phi}}\right)^p$ and $T_{e^{is\theta}r^d}=\left(T_{e^{i\theta}\widetilde{\psi}}\right)^s$. Now, if the hypotheses $(H1)$ and $(H2)$ are satisfied, then \cite[Remark 2]{lr} implies that for every vector $z^k$ of the orthogonal basis of $L^2_a(\mathbb{D})$, we have
\begin{eqnarray}
	T_{e^{im\theta}f}T_{e^{ip\theta}r^n}(z^k)&=&T_{e^{ip\theta}r^n}T_{e^{im\theta}f}(z^k),\\
	T_{e^{il\theta}g}T_{e^{is\theta}r^d}(z^k)&=&T_{e^{is\theta}r^d}T_{e^{il\theta}g}(z^k),\\
	\left(T_{e^{im\theta}f}T_{e^{is\theta}r^d}+T_{e^{il\theta}g}T_{e^{ip\theta}r^n}\right)(z^k)
	&=&\left(T_{e^{is\theta}r^d}T_{e^{im\theta}f}+T_{e^{ip\theta}r^n}T_{e^{il\theta}g}\right)(z^k).
\end{eqnarray}
From equations $(1)$ and $(2)$, we can conclude that $T_{e^{im\theta}f}$ and $T_{e^{il\theta}g}$ commute with $T_{e^{ip\theta}r^n}$ and  $T_{e^{is\theta}r^d}$, respectively. Therefore, \cite[Proposition 2 and Lemma 2]{lr} imply that
\begin{equation}\label{m}
T_{e^{im\theta}f}=c_1\left(T_{e^{i\theta}\widetilde{\phi}}\right)^m
\end{equation}
and
\begin{equation}\label{l}
	T_{e^{il\theta}g}=c_2\left(T_{e^{i\theta}\widetilde{\psi}}\right)^l
\end{equation}
for some constants $c_1$ and $c_2$. To avoid the trivial case where the operators are zero, we assume from now on that $c_1$ and $c_2$ are nonzero constants. Finally, considering the  commutator $[T,S]=TS-ST$ of two operators $T$ and $S$, Equation (3) can be written as
\begin{equation}\label{com}
	c_1\left[\left(T_{e^{i\theta}\widetilde{\phi}}\right)^m,T_{e^{is\theta}r^d}\right](z^k)=c_2\left[T_{e^{ip\theta}r^n},\left(T_{e^{i\theta}\widetilde{\psi}}\right)^l\right](z^k),\textrm{ for all }k\geq 0.
\end{equation}
{\noindent{{\bf{Important comments}}}
	\begin{itemize}
		\item[i)]  If $T_{e^{ip\theta}r^n}$ commutes with $T_{e^{is\theta}r^d}$, then \cite[Proposition 2 and Lemma 2]{lr} imply that $T_{e^{i\theta}\widetilde{\psi}}=cT_{e^{i\theta}\widetilde{\phi}}$ for some constant $c$. Moreover, \cite[Corollary 1]{lr} indicates that all four Toeplitz operators $T_{e^{ip\theta}r^n}$, $T_{e^{is\theta}r^d}$, $T_{e^{im\theta}f}$, and $T_{e^{il\theta}g}$ commute with each other. Consequently, they are all of the form a constant multiple of a single Toeplitz operator, which is $T_{e^{i\theta}\widetilde{\phi}}$. Thus, without loss of generality, we assume that $[T_{e^{ip\theta}r^n}, T_{e^{is\theta}r^d}]\neq 0$.
		\item[ii)] The case $p=s$ (and similarly, $l=m$) has been extensively studied and fully solved. For detailed discussions, see \cite{cr, l, lz}.
		\item[iii)] We will demonstrate that if $(H1)$ and $(H2)$ hold, then $m=p$, $l=s$, which implies that the constants $c_1$ and $c_2$ appearing in Equation (\ref{com}) are equal. Specifically, if $m=p$ (or equivalently, $l=s)$, then $(H2)$ ensures that $l=s$ (or equivalently,  $m=p)$. Moreover, Equation (\ref{m}) and Equation (\ref{l}) imply that $T_{e^{im\theta}f}=c_1T_{e^{ip\theta}r^n}$ and $T_{e^{il\theta}g}=c_2T_{e^{is\theta}r^d}$, where  $c_1$ and $c_2$ are constants. Consequently, Equation (\ref{com}) simplifies to
		$$c_1\left[T_{e^{ip\theta}r^n},T_{e^{is\theta}r^d}\right](z^k)=c_2\left[T_{e^{ip\theta}r^n},T_{e^{is\theta}r^d}\right](z^k),\textrm{ for all }k\geq 0.$$ Thus, $c_1=c_2$  assuming that $\left[T_{e^{ip\theta}r^n},T_{e^{is\theta}r^d}\right]\neq 0$.
		
	\end{itemize}
%\end{comments}
\section{Preliminaries and tools}
Radial functions in $L^1(\mathbb{D},dA)$ can be viewed as functions in $L^1([0,1),rdr)$. For a function $\phi\in L^1([0,1),rdr)$, we define its Mellin transform, denoted $\widehat{\phi}$, by  $$\widehat{\phi}(z)=\int_0^1\phi(r)r^{z-1}\,rdr.$$
It is well known that for  functions  $\phi\in L^1([0,1)rdr)$, the Mellin transform  $\widehat{\phi}$ is analytic on the right-half plane $\{z\in\mathbb{C}:\Re z>2\}$ and is continuous and bounded on $\{z\in\mathbb{C}:\Re z\geq2\}$.

The Mellin transform, which is related to the Laplace transform via the change of variable $r=e^{-t}$, is a valuable tool in studying quasihomogeneous Toeplitz operators. Indeed, quasihomogeneous Toeplitz operators act on the vectors of the orthogonal basis of $L^2_a(\mathbb{D})$ as shift operators with analytic weight, and this weight involves the Mellin transform of the symbol. We have the following lemma from \cite[Lemma 1, p.883]{lry}.
\begin{lemma} Let $\phi\in L^1([0,1), rdr)$ and let $p$ be a non-negative integer. Then for every integer $k\geq 0$, we have 
	$$T_{e^{ip\theta}\phi}(z^k)=2(k+p+1)\widehat{\phi}(2k+p+2)z^{k+p}.$$
\end{lemma}
Another important property of the Mellin transform is that $\widehat{\phi}$ is uniquely determined by its values on any set of integers satisfying the M\"{u}ntz-Sz\'asz (or Blaschke) condition. In our calculations, we
often determine $\phi$ by knowing its Mellin transform $\widehat{\phi}$ on an arithmetic sequences.  We have the following classical theorem \cite[p.102]{Rem}.
\begin{theorem}\label{thm3} 
Suppose that $f$ is a bounded analytic function on the right-half plane $\{z\in\mathbb{C}: \Re z>0\}$ that vanishes at a set of pairwise distinct points $d_1,d_2,\cdots,$ where
\begin{itemize}
\item[(i)]
$\inf \{\vert d_n\vert\}>0$, and
\item[(ii)]
$\sum\limits_{n\geq 1}^{} \Re(\frac{1}{d_n})=\infty$.
\end{itemize}
Then $f$ vanishes identically on $\{z\in\mathbb{C}:\Re z>0\}$.
\end{theorem}

The following lemma plays a key role in our  proof of the main result. Specifically, at a certain stage in the proof, we need to determine when the quotient of four Gamma functions is a rational function. We omit the proof of the lemma, which is a slight modification of \cite[Theorem 3, p.197-198]{cr}. 
\begin{lemma}\label{lem2}
Let $a, b, c, d$ be non-negative integers such that $a+b-c-d=\lambda$ and let $\delta \in\mathbb{N}$. Define the function $H$ to be
\begin{equation*}
H(z)=\frac{\Gamma\left(\frac{z}{2\delta}+\frac{a}{2\delta}\right)\Gamma\left(\frac{z}{2\delta}+\frac{b}{2\delta}\right)}{\Gamma\left(\frac{z}{2\delta}+\frac{c}{2\delta}\right)\Gamma\left(\frac{z}{2\delta}+\frac{d}{2\delta}\right)}.
\end{equation*}
Then, $H$ is a rational function if and only if $2\delta$ divides $\lambda$ and one of the numbers $a-c$ or $a-d$.
 
\end{lemma}
In \cite[Theorem 3, p.197-198]{cr}, the authors assume  $a+b-c-d=-1$ rather than $a+b-c-d=\lambda$ as in our version above. However, the proof remains exactly the same as stated in \cite[p.205]{cr}.  

\section{Main Results}
To ensure clarity and effectively persuade readers who may not be entirely familiar with the calculations in our proofs, we will begin by meticulously describing the case when 
$s=2p$. This approach aims to enhance understanding and reinforce confidence in the validity of our results. Moreover, the proof of the general case (Theorem 4) is fundamentally based on this key scenario.
\begin{theorem}\label{thm2}
Let $\phi(r)=r^n$ and $\psi(r)=r^d$ with $p<s$, $n,d\in\mathbb{N}$. Suppose there exist $m, l\in\mathbb{N}$ and nontrivial radial functions $f, g$ such that the hypotheses (H1) and (H2) of Theorem 1 are satisfied. If $s=2p$, then $m =p$, $l=s$ and
$$
T_{e^{im\theta}f}+T_{e^{is\theta}g}=c(T_{e^{ip\theta}r^n}+ T_{e^{is\theta}r^d}).
$$
for some constant $c$.
\end{theorem}
For simplicity, we will adopt the following notation in the proof: we will write $\equiv$ instead of $=$ when the quantity on the left side of the equation is equal to a constant multiple of the quantity on the right side.
\begin{proof}
Let $p, s, m, l, d, n\in\mathbb{N}^{*}$. Suppose that the hypotheses $(H1)$ and $(H2)$ are satisfied. In this case, Equation (\ref{com}) implies
\begin{align*}
	&\frac{2k + 2s + 2}{2k + s + d + 2} \prod_{j=0}^{m-1} \left(2k + 2s + 2j + 4\right) \widehat{\tilde{\phi}} \left(2k + 2s + 2j + 3\right) \\
	&- \frac{2k + 2m + 2s + 2}{2k + 2m + s + d + 2} \prod_{j=0}^{m-1} \left(2k + 2j + 4\right) \widehat{\tilde{\phi}} \left(2k + 2j + 3\right) \\
	&\equiv \frac{2k + 2l + 2p + 2}{2k + 2l + p + n + 2} \prod_{j=0}^{l-1} \left(2k + 2j + 4\right) \widehat{\tilde{\psi}} \left(2k + 2j + 3\right) \\
	&- \frac{2k + 2p + 2}{2k + p + n + 2} \prod_{j=0}^{l-1} \left(2k + 2p + 2j + 4\right) \widehat{\tilde{\psi}} \left(2k + 2p + 2j + 3\right).
\end{align*}

Set $z=2k+2$. By applying \cite[Theorem 14, p.1473]{l}, it follows that
\begin{align*}
	&\frac{(z + 2m + 2p)(z + p + n)(z + 2p + p + n) \Gamma\left(\frac{z + 2m}{2p}\right) \Gamma\left(\frac{z + p + n}{2p}\right)}{(z + s + d)(z + 2p)(z + 2m + p + n)(z + 2p + 2m + p + n) \Gamma\left(\frac{z}{2p}\right) \Gamma\left(\frac{z + 2m + p + n}{2p}\right)} \\
	&-\frac{1}{z + 2m + s + d} \cdot \frac{\Gamma\left(\frac{z + 2m}{2p}\right) \Gamma\left(\frac{z + p + n}{2p}\right)}{\Gamma\left(\frac{z}{2p}\right) \Gamma\left(\frac{z + 2m + p + n}{2p}\right)} \\
	&\equiv \frac{(z + 2l) \Gamma\left(\frac{z + 2l}{2s}\right) \Gamma\left(\frac{z + s + d}{2s}\right)}{(z + 2m) \Gamma\left(\frac{z}{2s}\right) \Gamma\left(\frac{z + 2l + s + d}{2s}\right)} \cdot \frac{1}{z + 2l + p + n} \\
	&- \frac{z \Gamma\left(\frac{z + 2m}{2s}\right) \Gamma\left(\frac{z + 2p + s + d}{2s}\right)}{(z + p + n)(z + 2m + s + d) \Gamma\left(\frac{z + 2p}{2s}\right) \Gamma\left(\frac{z + 2m + s + d}{2s}\right)}.
\end{align*}

This above equation is equivalent to
\begin{equation}\label{eq24}
	R_1(z) \cdot \frac{\Gamma\left(\frac{z + 2m}{2p}\right) \Gamma\left(\frac{z + p + n}{2p}\right)}{\Gamma\left(\frac{z}{2p}\right) \Gamma\left(\frac{z + 2m + p + n}{2p}\right)}
	\equiv R_2(z) \cdot \frac{\Gamma\left(\frac{z + 2l}{2s}\right) \Gamma\left(\frac{z + s + d}{2s}\right)}{\Gamma\left(\frac{z}{2s}\right) \Gamma\left(\frac{z + 2l + s + d}{2s}\right)}
	- R_3(z) \cdot \frac{\Gamma\left(\frac{z + 2m}{2s}\right) \Gamma\left(\frac{z + 2p + s + d}{2s}\right)}{\Gamma\left(\frac{z + 2p}{2s}\right) \Gamma\left(\frac{z + 2m + s + d}{2s}\right)},
\end{equation}

where $R_1, R_2$ and $R_3$ are rational function in $z$.
At this stage, we need to consider the case "$l<s$ and $m<p$" separately from the case "$l>s$ and $m>p$".
$${\framebox{{\bf{Case I}}. $l< s$ and $m< p$.}}$$
Our aim here is to show that the function $\dfrac{\Gamma\left(\frac{z + 2m}{2p}\right) \Gamma\left(\frac{z + p + n}{2p}\right)}{\Gamma\left(\frac{z}{2p}\right) \Gamma\left(\frac{z + 2m + p + n}{2p}\right)}
$ is rational. We will proceed by contradiction.  Assume, for the sake of contradiction, that this function 
is not rational. We consider the following two subcases:
\begin{itemize}
	\item[(1)]
	Suppose the functions $$\Gamma\left(\frac{z+2l}{2s}\right)\Gamma\left(\frac{z+s+d}{2s}\right)\textrm{ and } \Gamma\left(\frac{z+2m}{2s}\right)\Gamma\left(\frac{z+2p+s+d}{2s}\right)$$ have the same poles. According to Lemma \ref{lem2}, this implies that $2s$ divides $2m-s-d$. Hence,  there exists $N\in\mathbb{N}$ such that $s+d=2m+2sN$. Now, let $A$ be a sufficiently large integer. For such an $A$,  the number $-2p(2A+1)$ is not a pole for the functions $R_1$, $R_2$, or $R_3$. However, it
	is a zero of both the left side of Equation (\ref{eq24}) and the function $\dfrac{\Gamma\left(\frac{z+2m}{2s}\right)\Gamma\left(\frac{z+2p+s+d}{2s}\right)}{\Gamma\left(\frac{z+2p}{2s}\right)\Gamma\left(\frac{z+2m+s+d}{2s}\right)}$. Consequently, it is also a zero of the function $\dfrac{\Gamma\left(\frac{z+2l}{2s}\right)\Gamma\left(\frac{z+s+d}{2s}\right)}{\Gamma\left(\frac{z}{2s}\right)\Gamma\left(\frac{z+2l+s+d}{2s}\right)}$. This implies that there exists $B\in\mathbb{N}$ such that $-2p(2A+1)=-2sB-2l-s-d$, which leads to $2s(A-B)=2m+s+d$. We deduce that $2s$ divides $2m+s+d$, and consequently, $2s$ divides $4m$. This implies that $p$ divides $m$, which is a contradiction because $m<p$.
	\item[(2)]
	Suppose that the function $\dfrac{\Gamma\left(\frac{z+2l}{2s}\right)\Gamma\left(\frac{z+s+d}{2s}\right)}{\Gamma\left(\frac{z+2m}{2s}\right)\Gamma\left(\frac{z+2p+s+d}{2s}\right)}$ is not rational. Dividing both sides of Equation (\ref{eq24}) by $\Gamma\left(\frac{z+2m}{2s}\right)\Gamma\left(\frac{z+2p+s+d}{2s}\right)$ yields
	\begin{align*}
		&	R_1(z) \cdot \frac{\Gamma\left(\frac{z + 2m}{2p}\right) \Gamma\left(\frac{z + p + n}{2p}\right)}{\Gamma\left(\frac{z}{2p}\right) \Gamma\left(\frac{z + 2m + p + n}{2p}\right)\Gamma\left(\frac{z+2m}{2s}\right)\Gamma\left(\frac{z+2p+s+d}{2s}\right)} \\
		&\equiv R_2(z) \cdot \frac{\Gamma\left(\frac{z + 2l}{2s}\right) \Gamma\left(\frac{z + s + d}{2s}\right)}{\Gamma\left(\frac{z}{2s}\right) \Gamma\left(\frac{z + 2l + s + d}{2s}\right)\Gamma\left(\frac{z+2m}{2s}\right)\Gamma\left(\frac{z+2p+s+d}{2s}\right)}	-   \frac{R_3(z)}{\Gamma\left(\frac{z + 2p}{2s}\right) \Gamma\left(\frac{z + 2m + s + d}{2s}\right)},
	\end{align*}
	So, the following functions must be rational:
	$$
	\frac{\Gamma\left(\frac{z + 2m}{2p}\right) \Gamma\left(\frac{z + p + n}{2p}\right)}{\Gamma\left(\frac{z+2m}{2s}\right)\Gamma\left(\frac{z+2p+s+d}{2s}\right)}\quad\mbox{and}\quad
	\frac{\Gamma\left(\frac{z + 2l}{2s}\right) \Gamma\left(\frac{z + s + d}{2s}\right)}{\Gamma\left(\frac{z}{2s}\right)\Gamma\left(\frac{z+2l+s+d}{2s}\right)}.
	$$
	Applying Lemma \ref{lem2}, we find that $2s$ must divide $s+d$ and hence,
	the right side of Equation (\ref{eq24}) is rational. Therefore, it follows that the function $\dfrac{\Gamma\left(\frac{z+2m}{2p}\right)\Gamma\left(\frac{z+p+n}{2p}\right)}{\Gamma\left(\frac{z}{2p}\right)\Gamma\left(\frac{z+2m+p+n}{2p}\right)}$ is rational as well. This leads to a contradiction.
\end{itemize} 

Since we proved that the function $\dfrac{\Gamma\left(\frac{z + 2m}{2p}\right) \Gamma\left(\frac{z + p + n}{2p}\right)}{\Gamma\left(\frac{z}{2p}\right) \Gamma\left(\frac{z + 2m + p + n}{2p}\right)}
$ is rational,  Lemma \ref{lem2} implies that $2p$ must divide $p+n$. Consequently, $n$ must be an odd number multiple of $p$. 
Similarly, from Equation (\ref{eq24}), the functions $$\dfrac{\Gamma\left(\frac{z + 2l}{2s}\right) \Gamma\left(\frac{z + s + d}{2s}\right)}{\Gamma\left(\frac{z}{2s}\right) \Gamma\left(\frac{z + 2l + s + d}{2s}\right)}
 \textrm{ and }\dfrac{\Gamma\left(\frac{z + 2m}{2s}\right) \Gamma\left(\frac{z + 2p + s + d}{2s}\right)}{\Gamma\left(\frac{z + 2p}{2s}\right) \Gamma\left(\frac{z + 2m + s + d}{2s}\right)}
$$ must also be rational. Applying Lemma \ref{lem2} again, we find that $2s$ must divide $s+d$. Therefore, $d$ must be an odd number multiple of $s$. Thus, Theorem \ref{thm4} completes the proof. 

\newpage
$${\framebox{{\bf{Case II}}. $l> s$ and $m> p$.}}$$
Observe that if $p$ divides $m$, then the function $\dfrac{\Gamma\left(\frac{z+2m}{2p}\right)\Gamma\left(\frac{z+p+n}{2p}\right)}{\Gamma\left(\frac{z}{2p}\right)\Gamma\left(\frac{z+2m+p+n}{2p}\right)}$ is rational. In this case, we have the following two possibilities:
\begin{itemize}
\item[(1)]
If $m=2Np$, then the function $\dfrac{\Gamma\left(\frac{z+2l}{2s}\right)\Gamma\left(\frac{z+s+d}{2s}\right)}{\Gamma\left(\frac{z+2m}{2s}\right)\Gamma\left(\frac{z+2p+s+d}{2s}\right)}$ is rational. According to Lemma \ref{lem2}, $2s$ must divide $2m-s-d$. Since $2s$ divides $2m$, it follows that $2s$ must also divide $s+d$. Thus, there exists an integer $M\in\mathbb{N}$ such that   $2sM=s+d$. This implies that  $d$ is an odd number multiple of $s$.  Next, Equation (\ref{eq24}) becomes
\begin{equation*}
	R_1(z) \prod_{i=0}^{2N-1} \frac{z + 2pi}{z + p + n + 2pi} \equiv R_2(z) \prod_{i=0}^{M-1} \frac{z + 2si}{z + 2l + 2si} - R_3(z) \prod_{i=0}^{M-1} \frac{z + 2p + 2si}{z + 2m + 2si}.
\end{equation*}
The point $-p-n-2p(2N-1)$ is a pole of the right side of the equation above. This pole can only arise from  $-2l-p-n$ or $-p-n$, as all the others poles are multiples of $2p$. For the equation to hold, this pole should be canceled. Thus, $2p$ divides $p+n$, which implies that $n$ is an odd number multiple of $p$. Hence, Theorem 1 completes the proof.
\item[(2)]
Assume $m=(2N+1)p$. Then both $$\dfrac{\Gamma\left(\frac{z+2l}{2s}\right)\Gamma\left(\frac{z+s+d}{2s}\right)}{\Gamma\left(\frac{z}{2s}\right)\Gamma\left(\frac{z+2l+s+d}{2s}\right)} \textrm{ and }\dfrac{\Gamma\left(\frac{z+2m}{2s}\right)\Gamma\left(\frac{z+2p+s+d}{2s}\right)}{\Gamma\left(\frac{z+2p}{2s}\right)\Gamma\left(\frac{z+2m+s+d}{2s}\right)}$$ are rational. This leads to
\begin{equation}\label{eq25}
	R_1(z) \prod_{i=0}^{2N} \frac{z + 2pi}{z + p + n + 2pi} \equiv R_2(z) \prod_{i=0}^{N} \frac{z + 2si}{z + s + d + 2si} - R_3(z) \prod_{i=0}^{N-1} \frac{z + 2p + 2si}{z + 2p + s + d + 2si}.
\end{equation}
We shall argue by contradiction. Suppose $2s$ does not divide $s+d$. Then the two products on the right side of the previous equation cannot be simplified or canceled. Consider
the point $-4pN+4p$. This point is a zero of the function $R_3(z)\prod\limits_{i=0}^{N-1}\frac{(z+2p+2si)}{(z+2p+s+d+2si)}$. Hence, there must exist $\exists i\in\{0,\cdots,N-1\}$ such that $-4pN+4p=-2p-2si$, which leads to a contradiction because we are assuming that $s=2p$. Thus, our assumption that  $2s$ does not divide $s+d$ must be incorrect. Therefore, $2s$ does divide $s+d$, which implies that $d$ is an odd number multiple of $s$.\\
Now, let assume that $2p$ does not divide $p+n$. The point $-p-n-2p(2N)$ would then be a pole of the left side of Equation (\ref{eq25}), but all  poles on the right side are multiples of $2p$, leading is a contradiction. Therefore, $2p$ must divide $p+n$, which implies that $n$ is an odd number multiple of $p$. Hence, Theorem 1 completes the proof.
\end{itemize}
The proof of the case where $\frac{m}{p}\notin\mathbb{N}$ is analogous to \bf{Case I}. 

\end{proof}
For the following result, we will relax the condition previously required in Theorem 3, specifically $s=2p$.
\begin{theorem}\label{thm1}
Let $\phi(r)=r^n$ and $\psi(r)=r^d$ with $p<s$ and $n,d\in\mathbb{N}$. Suppose there exist $m, l\in\mathbb{N}$ and nontrivial radial functions $f, g$ such that the hypotheses (H1) and (H2) are satisfied. Then $m =p$, $l=s$ and
$$
T_{e^{im\theta}f}+T_{e^{is\theta}g}=c(T_{e^{ip\theta}r^n}+ T_{e^{is\theta}r^d}),
$$
for some constant $c$.
\end{theorem}

\begin{proof}
By setting $z=2k+2$ in Equation (6) and applying \cite[Theorem 14, p.1473]{l}, we obtain:
\begin{align*}
	&\frac{z + 2s}{z + s + d} \cdot \frac{(z + 2s + 2m) \Gamma\left(\frac{z + 2s + 2m}{2p}\right) \Gamma\left(\frac{z + 2s + p + n}{2p}\right)}{(z + 2s) \Gamma\left(\frac{z + 2s}{2p}\right) \Gamma\left(\frac{z + 2s + 2m + p + n}{2p}\right)} \\
	&- \frac{z + 2m + 2s}{z + 2m + s + d} \cdot \frac{(z + 2m) \Gamma\left(\frac{z + 2m}{2p}\right) \Gamma\left(\frac{z + p + n}{2p}\right)}{z \Gamma\left(\frac{z}{2p}\right) \Gamma\left(\frac{z + 2m + p + n}{2p}\right)} \\
	&\equiv \frac{z + 2l + 2p}{z + 2l + p + n} \cdot \frac{(z + 2l) \Gamma\left(\frac{z + 2l}{2s}\right) \Gamma\left(\frac{z + s + d}{2s}\right)}{z \Gamma\left(\frac{z}{2s}\right) \Gamma\left(\frac{z + 2l + s + d}{2s}\right)} \\
	&- \frac{z + 2p}{z + p + n} \cdot \frac{(z + 2p + 2l) \Gamma\left(\frac{z + 2p + 2l}{2s}\right) \Gamma\left(\frac{z + 2p + s + d}{2s}\right)}{(z + 2p) \Gamma\left(\frac{z + 2p}{2s}\right) \Gamma\left(\frac{z + 2p + 2l + s + d}{2s}\right)}.
\end{align*}

By (H2), there exist $\alpha, \beta\in\mathbb{N}$ such that $m+\beta=l$ and $p+\alpha=s$. Since $l+p=m+s$, we have $m+\beta+p=m+p+\alpha$. Therefore, $\beta=\alpha$. Let $\alpha\in\mathbb{N}$ be such that $p+\alpha=s, m+\alpha=l$. Then Equation (\ref{com}) becomes
\begin{align*}
	&\frac{(z+2s)}{(z+s+d)} \cdot \frac{(z+2s+2m)(z+2\alpha+2m)(z+2\alpha+p+n) \Gamma\left(\frac{z+2\alpha+2m}{2p}\right) \Gamma\left(\frac{z+2\alpha+p+n}{2p}\right)}{(z+2s)(z+2\alpha)(z+2\alpha+2m+p+n) \Gamma\left(\frac{z+2\alpha}{2p}\right) \Gamma\left(\frac{z+2\alpha+2m+p+n}{2p}\right)} \\
	&- \frac{(z+2m+2s)}{(z+2m+s+d)} \cdot \frac{(z+2m) \Gamma\left(\frac{z+2m}{2p}\right) \Gamma\left(\frac{z+p+n}{2p}\right)}{z \Gamma\left(\frac{z}{2p}\right) \Gamma\left(\frac{z+2m+p+n}{2p}\right)} \\
	&\equiv \frac{(z+2l+2p)}{(z+2l+p+n)} \cdot \frac{(z+2l) \Gamma\left(\frac{z+2l}{2s}\right) \Gamma\left(\frac{z+s+d}{2s}\right)}{z \Gamma\left(\frac{z}{2s}\right) \Gamma\left(\frac{z+2l+s+d}{2s}\right)} \\
	&- \frac{(z+2p)}{(z+p+n)} \cdot \frac{(z+2p+2l)(z+2m) \Gamma\left(\frac{z+2m}{2s}\right) \Gamma\left(\frac{z+2p+s+d}{2s}\right)}{(z+2p)(z+2m+s+d) \Gamma\left(\frac{z+2p}{2s}\right) \Gamma\left(\frac{z+2m+s+d}{2s}\right)}.
\end{align*}
Thus, we obtain
\begin{align*}
	&\frac{(z+2\alpha+p+n)(z+2m+s+d)(z+2l)}{(z+s+d)(z+2l+p+n)(z+2\alpha)} \Bigg[ \frac{\Gamma\left(\frac{z+2\alpha+2m}{2p}\right) \Gamma\left(\frac{z+2\alpha+p+n}{2p}\right)}{\Gamma\left(\frac{z+2\alpha}{2p}\right) \Gamma\left(\frac{z+2\alpha+2m+p+n}{2p}\right)} \\
	&\quad - \frac{(z+s+d)(z+2\alpha) \Gamma\left(\frac{z+2l}{2s}\right) \Gamma\left(\frac{z+s+d}{2s}\right)}{(z+2\alpha+p+n)(z) \Gamma\left(\frac{z}{2s}\right) \Gamma\left(\frac{z+2l+s+d}{2s}\right)} \Bigg] \\
	&\equiv \frac{z+2m}{z} \left[ \frac{\Gamma\left(\frac{z+2m}{2p}\right) \Gamma\left(\frac{z+p+n}{2p}\right)}{\Gamma\left(\frac{z}{2p}\right) \Gamma\left(\frac{z+2m+p+n}{2p}\right)} - \frac{z \Gamma\left(\frac{z+2m}{2s}\right) \Gamma\left(\frac{z+2p+s+d}{2s}\right)}{(z+p+n) \Gamma\left(\frac{z+2p}{2s}\right) \Gamma\left(\frac{z+2m+s+d}{2s}\right)} \right].
\end{align*}
This simplifies to
\begin{equation}\label{eq21}
H(z)F(z+2\alpha)\equiv F(z),
\end{equation}
where
\begin{equation*}
	F(z) = \frac{z+2m}{z} \left[ \frac{\Gamma\left(\frac{z+2m}{2p}\right) \Gamma\left(\frac{z+p+n}{2p}\right)}{\Gamma\left(\frac{z}{2p}\right) \Gamma\left(\frac{z+2m+p+n}{2p}\right)} - \frac{z \Gamma\left(\frac{z+2m}{2s}\right) \Gamma\left(\frac{z+2p+s+d}{2s}\right)}{(z+p+n) \Gamma\left(\frac{z+2p}{2s}\right) \Gamma\left(\frac{z+2m+s+d}{2s}\right)} \right].
\end{equation*}

and 
\begin{equation*}
H(z) = \dfrac{(z+2\alpha+p+n)(z+2m+s+d)}{(z+2l+p+n)(z+s+d)}.
\end{equation*}
It is easy to see that the function $F$ has infinitely many poles. Specifically, $F$ has poles at the points $-2sA-2m, -2sB-2p-s-d, -2pC-2m, -2pD-p-n$ for sufficiently large integers $A, B, C, D\in\mathbb{N}$. These integers  are chosen large enough to ensure that these poles are not canceled out by the zeros of the function $H(z)$. These poles arise from the terms $\Gamma\left(\frac{z+2m}{2s}\right)\Gamma\left(\frac{z+2p+s+d}{2s}\right)$ and $\Gamma\left(\frac{z+2m}{2p}\right)\Gamma\left(\frac{z+p+n}{2p}\right)$. Additionally, Equation (\ref{eq21}) indicates that these points are also poles of the function $F(z+2\alpha)$. We then consider the following two situations:
\begin{itemize}
\item[\textbf{Situation 1:}]
The poles of $\Gamma\left(\frac{z+2m}{2s}\right)\Gamma\left(\frac{z+2p+s+d}{2s}\right)$ arise from $\Gamma\left(\frac{z+2\alpha+2m}{2s}\right)\Gamma\left(\frac{z+s+d}{2s}\right)$. Similarly, the poles of $\Gamma\left(\frac{z+2m}{2p}\right)\Gamma\left(\frac{z+p+n}{2p}\right)$ arise from $\Gamma\left(\frac{z+2\alpha+2m}{2p}\right)\Gamma\left(\frac{z+2\alpha+p+n}{2p}\right)$. This implies that the functions $$\dfrac{\Gamma\left(\frac{z+2m}{2s}\right)\Gamma\left(\frac{z+2p+s+d}{2s}\right)}{\Gamma\left(\frac{z+2\alpha+2m}{2s}\right)\Gamma\left(\frac{z+s+d}{2s}\right)} \textrm{ and } \dfrac{\Gamma\left(\frac{z+2m}{2p}\right)\Gamma\left(\frac{z+p+n}{2p}\right)}{\Gamma\left(\frac{z+2\alpha+2m}{2p}\right)\Gamma\left(\frac{z+2\alpha+p+n}{2p}\right)}$$ are rational functions. By applying Lemma \ref{lem2}, we conclude that  $s$ must divide $p-\alpha$ and $p$ must divide $2\alpha$. Hence, $p=\alpha$. Since $s=p+\alpha$, it follows that $s=2p$. Therefore, Theorem \ref{thm2} completes the proof. 
\item[\textbf{Situation 2:}]
The poles of $\Gamma\left(\frac{z+2m}{2s}\right)\Gamma\left(\frac{z+2p+s+d}{2s}\right)$ arise from $\Gamma\left(\frac{z+2\alpha+2m}{2p}\right)\Gamma\left(\frac{z+2\alpha+p+n}{2p}\right)$. Similarly, the poles of $\Gamma\left(\frac{z+2m}{2p}\right)\Gamma\left(\frac{z+p+n}{2p}\right)$ arise from $\Gamma\left(\frac{z+2\alpha+2m}{2s}\right)\Gamma\left(\frac{z+s+d}{2s}\right)$. This implies  that the functions $$\dfrac{\Gamma\left(\frac{z+2m}{2s}\right)\Gamma\left(\frac{z+2p+s+d}{2s}\right)}{\Gamma\left(\frac{z+2\alpha+2m}{2p}\right)\Gamma\left(\frac{z+2\alpha+p+n}{2p}\right)} \textrm{ and } \dfrac{\Gamma\left(\frac{z+2m}{2p}\right)\Gamma\left(\frac{z+p+n}{2p}\right)}{\Gamma\left(\frac{z+2\alpha+2m}{2s}\right)\Gamma\left(\frac{z+s+d}{2s}\right)}$$ are rational functions. Furthermore, the function $\dfrac{\Gamma\left(\frac{z+2m}{2p}\right)\Gamma\left(\frac{z+p+n}{2p}\right)}{\Gamma\left(\frac{z+2m}{2s}\right)\Gamma\left(\frac{z+2p+s+d}{2s}\right)}$ is also rational. Consequently, when  multiplying Equation (9) with the ratio $\dfrac{\Gamma\left(\frac{z+2\alpha+2m}{2s}\right)\Gamma\left(\frac{z+s+d}{2s}\right)}{\Gamma\left(\frac{z+2m}{2p}\right)\Gamma\left(\frac{z+p+n}{2p}\right)}$, we obtain the function  $$\dfrac{\Gamma\left(\frac{z+2\alpha+2m}{2s}\right)\Gamma\left(\frac{z+s+d}{2s}\right)}{\Gamma\left(\frac{z+2m}{2s}\right)\Gamma\left(\frac{z+2p+s+d}{2s}\right)},$$ which is also rational. According to Lemma \ref{lem2}, this implies that $s$ divides $\alpha-p$. If $\alpha>p$, then $s=p+\alpha\leq \alpha -p$ which is not possible. Conversely, if $\alpha< p$, then $s=p+\alpha\leq p-\alpha$ which is also not possible. Thus, $p=\alpha$,  which implies $s=2p$. Hence, Theorem \ref{thm2} completes the proof.
\end{itemize}

\end{proof}
%\newpage
\section*{Acknowledgment} The authors would like to express their gratitude to the reviewers for their invaluable comments, which significantly enhanced the presentation of this paper.


\begin{thebibliography}{1}
\bibitem{sl} Al Sabi, Hashem and Louhichi, Issam, \textit{On the commutativity of Toeplitz operators with harmonic symbols}, Oper. Matrices 12 (2018), no. 4, 1159-1176.
\bibitem{al} Aqel, Khitam and Louhichi, Issam, \textit{On the commutativity of sums of Toeplitz operators on the Bergman space}, Grad. J. Math. 2 (2017), no. 2, 51-58.
\bibitem{acr} Axler, Sheldon; \u{C}u\u{c}kovi\'{c}, \u{Z}eljko; and Rao, N. V., \textit{Commutants of analytic Toeplitz operators on the Bergman space}, Proc. Amer. Math. Soc. 128 (2000), no. 7, 1951-1953.
\bibitem{cr} \u{Z}. \u{C}u\u{c}kovi\'{c} and N. V. Rao, \textit{Mellin transform, monomial symbols, and commuting Toeplitz operators}, J. Funct. Anal. 154 (1998), no. 1, 195-214. 
\bibitem{h} Hedenmalm, Haakan; Korenblum, Boris; and Zhu, Kehe, \textit{Theory of Bergman spaces}, Graduate Texts in Mathematics, 199. Springer-Verlag, New York, 2000. x+286 pp. ISBN: 0-387-98791-6
\bibitem{lt} Le, Trieu and Tikaradze, Akaki, \textit{Commutants of Toeplitz operators with harmonic symbols}, New York J. Math. 23 (2017), 1723-1731.
\bibitem{l} Louhichi, Issam, \textit{Powers and roots of Toeplitz operators}, Proc. Amer. Math. Soc. 135 (2007), no. 5, 1465-1475.
\bibitem{lrp} Louhichi, Issam and Rao, N. V., \textit{Roots of Toeplitz operators on the Bergman space}, Pacific J. Math. 252 (2011), no. 1, 127-144.
\bibitem{lr} Louhichi, I. and Rao, N. V., \textit{Bicommutants of Toeplitz operators}, Arch. Math. (Basel) 91 (2008), no. 3, 256-264.
\bibitem{lry} Louhichi, Issam; Rao, Nagisetti V.; and Yousef, Abdel, \textit{Two questions on products of Toeplitz operators on the Bergman space}, Complex Anal. Oper. Theory 3 (2009), no. 4, 881-889.
\bibitem{lsz} I. Louhichi, E. Strouse, and L. Zakariasy, \textit{Products of Toeplitz operators on the Bergman space}, Integral Equations and Operator Theory \textbf{54} No. \textbf{4}, p 525-539, 2006.
\bibitem{lz} Louhichi, Issam and Zakariasy, Lova, \textit{Quasihomogeneous Toeplitz operators on the harmonic Bergman space}, Arch. Math. (Basel) 98 (2012), no. 1, 49-60.
\bibitem{Rem} R. Remmert, \textit{Classical Topics in complex function theory}, Graduate Texts in Mathematics \textbf{172}, Springer, New York, 1998. 
\bibitem{rv} Rozenblum, G., Vasilevski, N., \textit{Commutative Algebras of Toeplitz Operators on the Bergman Space Revisited: Spectral Theorem Approach}, Integr. Equ. Oper. Theory 94, 27 (2022). 
\bibitem{sz} Strouse, Elizabeth and Zakariasy, Lova, \textit{Quasihomogeneous Toeplitz operators}, Harmonic analysis, function theory, operator theory, and their applications, 251-260, Theta Ser. Adv. Math., 19, Theta, Bucharest, 2017.
\bibitem{y} Yousef, A. and Al-Naimi, R., \textit{On Toeplitz Operators with Biharmonic Symbols}, Bull. Malays. Math. Sci. Soc. 43, 1647–1659 (2020). https://doi.org/10.1007/s40840-019-00763-3




\end{thebibliography}
\end{document}